\documentclass[reqno,11pt]{amsart}
\usepackage[utf8]{inputenc}
\usepackage[T1,T5]{fontenc}
\usepackage{amsmath,amsfonts,amssymb,amsthm}
\usepackage{appendix}
\usepackage[left=2cm,right=2cm,top=2cm,bottom=2.5cm]{geometry}
\usepackage{stmaryrd} 
\usepackage{mathrsfs}
\usepackage{subfiles} 
\usepackage{yhmath} 
\usepackage{cases} 
\usepackage{mathtools,color,url} 
\usepackage{abraces}
\usepackage{bbm} 
\usepackage{tikz}

\newtheorem{theo}{Theorem}
\newtheorem{lemm}[theo]{Lemma}
\newtheorem{coro}[theo]{Corollary}
\newtheorem{defi}[theo]{Definition}

\newtheorem{claim}[theo]{Claim}
\newtheorem*{stat}{Statement}
\newtheorem{rema}[theo]{Remark}

\newcommand\R{{\ensuremath {\mathbb R} }}

\newcommand\N{{\ensuremath {\mathbb N} }}

\newcommand{\sign}{\text{sign}} 

\mathtoolsset{showonlyrefs}  

\DeclareMathOperator{\sech}{sech}

\title{Decay of solitary waves of fractional Korteweg-de Vries type equations}
\author[A.Eychenne]{Arnaud Eychenne}
\address{Department of Mathematics, University of Bergen, Allégaten 41
Realfagbygget
5007 Bergen, Norway}
\email{arnaud.eychenne.waxweiler@gmail.com}
\author[F.Valet]{Frédéric Valet}
\address{Department of Mathematics, University of Bergen, Allégaten 41
Realfagbygget
5007 Bergen, Norway}
\email{frederic.valet@uib.no}
\date{\today}


\begin{document}
\maketitle
\begin{abstract}
We study the solitary waves of fractional Korteweg-de Vries type equations, that are related to the $1$-dimensional semi-linear fractional equations:
\begin{align*}
    \vert D \vert^\alpha u + u -f(u)=0,
\end{align*}
with $\alpha\in (0,2)$, a prescribed coefficient $p^*(\alpha)$, and a non-linearity $f(u)=\vert u \vert^{p-1}u$ for $p\in(1,p^*(\alpha))$, or $f(u)=u^p$ with an integer $p\in[2;p^*(\alpha))$. Asymptotic developments of order $1$ at infinity of solutions are given, as well as second order developments for positive solutions, in terms of the coefficient of dispersion $\alpha$ and of the non-linearity $p$. The main tools are the kernel formulation introduced by Bona and Li, and an accurate description of the kernel by complex analysis theory.
\end{abstract}

\section{Introduction}

\subsection{Motivation}

Consider the class of semi-linear fractional equations:
\begin{align}\label{defi:eq_elliptic}
    \vert D \vert^\alpha u + u -f(u)=0, \quad u:\R\rightarrow\R, \quad 0 < \alpha < 2, 
\end{align}
with the notation $\vert D \vert^{\alpha}$  standing for the Riesz potential of order $-\alpha$ given by the Fourier multiplier:
\begin{align*}
    \mathcal{F}(\vert D \vert^{\alpha} u ) := \vert \xi \vert^\alpha \mathcal{F}(u),
\end{align*}
and with a non-linearity $f(u)=\vert u \vert^{p-1} u$ or $f(u)= u^p$, $p$ integer, and $1< p< p^*(\alpha)$ where:
\begin{equation*}
    p^*(\alpha):=\left\{
    \begin{aligned}
    &\frac{2\alpha}{1-\alpha}+1 \quad & \text{if } 0 < \alpha \leq 1, \\
    &+\infty \quad & \text{if } 1\leq \alpha < 2.
    \end{aligned}
    \right.
\end{equation*}
$p^*(\alpha)$ is defined such that the "equation is $H^{\frac{\alpha}{2}}(\R)$-subcritical", where $H^{s}(\R)$ stands for the Sobolev spaces.

Those fractional elliptic equations appear naturally when studying solitary waves of the following equations:
\begin{itemize}
    \item the fractional generalized Korteweg-de Vries equation \cite{SV96,SW21}:
\begin{align}\label{fgKdV}\tag{fgKdV}
    \partial_t u +\partial_x \left( -\vert D \vert^\alpha u + f(u) \right)=0,
\end{align}
    \item the fractional generalized Benjamin-Bona-Mahony equation \cite{EEE15,EEE16,OBM20}:
\begin{align}\label{fgBBM}\tag{fgBBM}
    \partial_t u +\partial_x u + \partial_x (f(u))-\vert D \vert^\alpha \partial_t u- \vert D \vert^\alpha \partial_x u  =0
\end{align}
    \item the fractional generalized nonlinear Schr\"{o}dinger equation \cite{Las02}:
    \begin{align}\label{fgNLS}\tag{fgNLS}
        i\partial_tu + \vert D \vert^\alpha u +f(u)=0. 
    \end{align}
\end{itemize}

The particular case of $f(u)=u^p$ with an integer $p$ is particularly relevant for \eqref{fgKdV} and \eqref{fgBBM}, whereas the non-linearity $f(u)=\vert u \vert^{p-1}u $ naturally appears when studying the Schr\"{o}dinger type equations.

A solitary wave for \eqref{fgKdV} is a wave moving in one direction with a constant velocity $c$, keeping its form along the time and decaying at infinity, and can thus be written as $u(t,x)=Q_c(x-ct)$. In the case $c=1$, the function $Q_1$ has to be a solution of \eqref{defi:eq_elliptic}.

Up to adequate change of variables, the solitary waves of \eqref{fgBBM} satisfy the same equation. The counterparts of \eqref{fgNLS} are stationary waves of the form $e^{i\omega t}Q_\omega(x)$, where the solution $Q_1$ with phase $\omega=1$ also satisfies \eqref{defi:eq_elliptic}.

Having a deeper understanding of the solutions of the semi-linear fractional equation \eqref{defi:eq_elliptic} is necessary to get more insights on the behaviour of solitary waves. It is of great interest to get more properties of those solutions, such as the existence and uniqueness of solutions, their regularity, the number of zeros or the asymptotic behaviour. The aim of this article to give the asymptotic development of the solutions to that semi-linear fractional equations and some subsequent properties.

\subsection{Survey of properties of the solutions.}

In the two specific cases $\alpha=2$ for general $p$, and $\alpha=1$ with $f(u)=u^2$, the solutions of the semi-linear fractional equation \eqref{defi:eq_elliptic} have been widely studied during the last fifty years. The main questions on those solutions are their existence, uniqueness, and some intrinsic properties.

Let us give a brief review of the previous results in the case $\alpha=2$, for algebraic non-linearities $p\in (1,p^*(\alpha))$. We consider the set of solutions $Q:\R^d\rightarrow \R$ with $\vert Q(x) \vert \underset{x\longrightarrow +\infty}{\rightarrow} 0$. For the dimension $d=1$, there exists a unique solution $Q(x)= \left(\frac{(p+1)}{2}\sech^{2}\left(\frac{p-1}{2} \right)\right)^{\frac{1}{p-1}}$ up to translations of the origin (see Theorem 5 of Berestycki-Lions \cite{BL83}). For higher dimensions $d\geq 2$, there exists a unique positive solution (see Weinstein \cite{weinstein1987existence} and Kwong \cite{Kwo89}), but also an infinite number of non-positive solutions as shown in Strauss \cite{S77},\cite{BL83} and Musso-Pacard-Wei \cite{MPW12}. For a more general local setting, we refer to the work of Adachi-Shibata-Watanabe \cite{ASW17} for uniqueness of ground states of quasi-linear elliptic equations.  All the previous solutions have an exponential decay.

Concerning the equation associated with the Benjamin-Ono equation (BO), corresponding to $\alpha=1$, $f(Q)=Q^2$ and $d=1$, a solution is explicit $Q(x)=4(1+x^2)^{-1}$, see Benjamin \cite{Ben67}. In fact, any solution to this problem is equal, up to translation, to this solution, as proved in Amick-Toland \cite{AT91Uniqueness} relying on former ideas of Benjamin  \cite{Ben67} (see also Albert \cite{Alb95} for an alternative proof). Notice that the solution is even, positive and has a polynomial decay at infinity.

In the generic case $\alpha\in(0,2)$, the existence of solutions relies on the existence of a minimizer for the functional $J^\alpha$ defined by:
\begin{align}\label{GagliardoNirenberg}
		\displaystyle J^{\alpha}(u)=\frac{\left(\displaystyle\int ||D|^{\frac{\alpha}{2}}u|^2\right)^{\frac{p+1}{2\alpha}}\left(\displaystyle\int |u|^2 \right)^{\frac{(p+1)}{2\alpha}(\alpha-1)+1}}{\displaystyle\int \vert u\vert ^{p+1}}.
\end{align} 
A minimizer $Q$ of the functional is called a ground state. The existence of a minimizer has been obtained by Weinstein \cite{weinstein1987existence} and Albert-Bona-Saut \cite{albert1997model}. The structure of the set of solutions of \eqref{defi:eq_elliptic} is complex, and it is not easy to know which elements compose this set. For example, the question of existence and uniqueness of solutions of \eqref{defi:eq_elliptic} not minimizing the functional $J^{\alpha}$ remains open. However, a breakthrough was achieved by Frank-Lenzmann \cite{FL13} by proving the uniqueness of ground states. Their proof relies on the non degeneracy of the linearized operator $L=|D|^{\alpha}+1-pQ^{p-1}$, in other words ker$(L)=$span$(Q')$. The understanding of the kernel of $L$ is based on a result by Caffarelli-Silvestre \cite{CS07} to express $|D|^{\alpha}$ as a Dirichlet-to-Neumann operator for a local problem on the upper half-plane. Furthermore in higher dimensions, Felmer-Quaas-Tan in \cite{FQT12} derived the existence and some properties of the ground states  and Frank-Lenzmann-Silvestre in \cite{FLS16} extended the uniqueness result. 

In the following theorem, we summarize the well-known properties of the ground states.

\begin{theo}[\cite{albert1997model,FL13,FLS16,weinstein1987existence}]\label{theo:resum_Q}
	Let $\alpha\in(0,2)$ and $p\in (1,p^{*}(\alpha))$. There exists $Q\in H^{s}(\R)$, for all $s\geq0$, such that  
	\begin{enumerate}
	\item (\textit{Existence}) The function $Q$ solves \eqref{defi:eq_elliptic} and  $Q=Q(|x|)>0$ is even, positive and strictly decreasing in $|x|$. Moreover, the function $Q$ is a minimizer of $J^{\alpha}$ in the sense that: 
		\begin{equation} 
		J^{\alpha}(Q)=\inf_{u \in H^{\frac{\alpha}2}(\mathbb R)}J^{\alpha}(u).
		\end{equation}	
	\item (\textit{Uniqueness}) The even ground state solution $Q=Q(|x|)>0$ of \eqref{defi:eq_elliptic} is unique. Furthermore, every optimizer $v\in H^{\frac{\alpha}{2}}(\R)$ for the Gagliardo-Nirenberg problem \eqref{GagliardoNirenberg} is of the form $v=\beta Q(\gamma(\cdot+y))$ with some $\beta\in\R, \beta\neq0, \gamma>0$ and $y\in \R$.
	\item (\textit{Decay}) The function $Q$ verifies the following decay estimate:
		\begin{align}\label{eq:decay}
			\frac{C_1}{\langle x \rangle^{\alpha+1}} \leq Q(x)\leq \frac{C_2}{\langle x \rangle^{\alpha+1}},
		\end{align} 
     for some $C_1,C_2>0$.
	\end{enumerate}
\end{theo}

As explained above, it is not clear whether other solutions of \eqref{defi:eq_elliptic} which are not minimizers of \eqref{GagliardoNirenberg} exist or not. The critical points of the functional $J^\alpha$ are called bound states and solve \eqref{defi:eq_elliptic}. In dimension $1$ for some elliptic equations, we expect the ground state to be the unique solution vanishing at infinity. For some elliptic equations in higher dimensions, other bound states different from the ground state exist: they are called excited states. For non-local equations such as \eqref{defi:eq_elliptic}, except for $\alpha=1$ and $f(u)=u^2$, the existence of a unique bound state in dimension $1$ and of multiple excited states in higher dimensions is still an open problem. In this context, we propose in this article to sharpen the asymptotic behaviour of the ground-state and to address the issue of the asymptotic behaviour of any solution of \eqref{defi:eq_elliptic} at $+\infty$.

\subsection{Main results} 

In this paper, we give several results on the asymptotic expansion of a solution $Q$ of \eqref{defi:eq_elliptic} and of its derivatives, and extend this development in the case of a non-linearity $f(u)=u^3$. These results extend the ones of Cappiello-Gramchev-Rodino \cite{CGR15}, where they proved that the solutions of a wider class of equations remain in some algebraic weighted spaces.

All the theorems stated in this article are given for $x>1$. However the proofs can be adapted to get the asymptotic developments at $-\infty$. 

To do so, we define the following explicit constants:
\begin{align}\label{defi:k1_k2}
    k_1:=\frac{\sin\left(\frac{\pi}{2}\alpha\right)}{\pi}\Gamma\left(\alpha+1 \right), \quad k_2:=-\displaystyle\frac{\sin\left(\pi\alpha \right)}{\pi} \Gamma\left(2\alpha+1\right) 
\end{align}
with $\Gamma$ the Euler $\Gamma$-function, and fir a given function $Q$:
\begin{align}\label{defi:a1_a2_a3}
    a_1:= k_1 \int \vert Q \vert^{p-1}Q(x) dx, \quad a_2:=k_2 \int \vert Q \vert^{p-1}Q(x) dx, \quad a_3 := \displaystyle \frac{(\alpha+1)(\alpha+2)}{2}k_1 \int x^2 \vert Q \vert^{p-1} Q(x) dx.
\end{align}

We recall the definition of a weak solution of \eqref{defi:eq_elliptic}:

\begin{defi}
    The function $u$ is called a weak solution of \eqref{defi:eq_elliptic} if for all $\varphi$ in the Schwartz space $\mathcal{S}(\R)$, we have that:
    \begin{align*}
        \int_\R \left(  u(x)|D|^{\alpha}\varphi(x) + u(x)\varphi(x) -f(u)(x)\varphi(x) \right)dx=0. 
    \end{align*}
\end{defi}

The next theorem states the main order terms of the development of the derivatives.

\begin{theo}\label{propo:notre_method_1}
Let $\alpha\in(0,2)$, $p\in(1,p^{*}(\alpha))$ and $Q$ be a weak solution of $|D|^{\alpha}Q+Q-|Q|^{p-1}Q=0$, satisfying:
\begin{align}\label{hypo:Q}
Q\in L^p(\mathbb{R}) \quad \text{and} \quad \exists l>0, \quad \vert x \vert^l Q(x) \in L^\infty(\mathbb{R}).
\end{align}
Then, $Q\in C^{0}(\R)$ and verifies:
\begin{equation*}
        \exists \beta>\alpha+1, \quad  Q(x) - \frac{a_1}{x^{\alpha+1}} =o_{+\infty} \left(x^{-\beta}\right),
\end{equation*}
    with $a_1$ dependent on $Q$ and defined in \eqref{defi:a1_a2_a3}.
    
Furthermore, if $\alpha>1$, then $Q\in C^{\lfloor p \rfloor+1}(\R)$ with $\lfloor p \rfloor$ the floor function of $p$, and verifies for $j\leq \lfloor p \rfloor $:
    \begin{align*}
        \exists \beta =\beta(j)>\alpha+1+j, \quad Q^{(j)}(x) -(-1)^{j}\frac{(\alpha+j)!}{\alpha!} \frac{a_1}{x^{\alpha+1+j}}  =o_{+\infty} \left(x^{-\beta}\right).
    \end{align*}
\end{theo}

\begin{rema} Some comments on the previous result are in order.
\begin{enumerate}
    \item Bona-Li in \cite{BL97} studied the decay of solutions of elliptic equations similar to \eqref{defi:eq_elliptic}, with a set of assumptions different from \eqref{hypo:Q}.  In the context of \eqref{defi:eq_elliptic}, if we ask for $\alpha>\frac{3}{2}$, $Q \in L^\infty(\R)$ and vanishing at infinity, Theorem 3.1.2 of \cite{BL97} implies the condition \eqref{hypo:Q}. Thus condition \eqref{hypo:Q} is coherent with their article.
    \item Note that if $Q$ is not positive the coefficient $a_1$ is potentially null.
    \item If $\alpha>1$,  our method does not give an asymptotic expansion of $Q^{(\lfloor p \rfloor+1)}$. Moreover, we do not know if $Q^{(\lfloor p \rfloor+2)}$ exists.
    \item If the coefficient of the non-linearity $p$ is an integer then $Q\in H^{\infty}(\R)$. The proof is given in Appendix \ref{appendix_A}.
\end{enumerate}
\end{rema}

In the next theorem, we assume the function $Q$ to be positive, so that Theorem 
\ref{theo:resum_Q} applies. We recall some important results of this theorem and give new asymptotic behaviours.

\begin{theo}\label{propo:notre_method_2}
Let $\alpha\in(0,2)$. Let $Q$ satisfying the assumptions \eqref{hypo:Q} of Theorem \ref{propo:notre_method_1}. Suppose also that $Q$ is positive. Then $Q\in H^{\infty}(\R)$, even (up to translation), decaying  and verifies that:
    \begin{align*}
        \forall j\in\N, \exists \beta=\beta(j) >\alpha+1+j, \quad Q^{(j)}(x) -(-1)^{j}\frac{(\alpha+j)!}{\alpha!} \frac{a_1}{x^{\alpha+1+j}} = o_{+\infty} \left(x^{-\beta}\right),
    \end{align*}
    and the next order asymptotic expansion holds, with a positive constant $\tilde{a}_1$:
\begin{equation*}
    \begin{aligned}
        & \text{Case } p<\frac{2\alpha+1}{\alpha+1}: \quad \exists \beta> p(\alpha+1), & Q(x) - \frac{a_1}{x^{\alpha+1}} - & \frac{\tilde{a}_1}{x^{p(\alpha+1)}} =o_{+\infty} \left(x^{-\beta}\right). \\
        & \text{Case } p=\frac{2\alpha+1}{\alpha+1}: \quad \exists \beta> 2\alpha+1, & Q(x)- \frac{a_1}{x^{\alpha+1}} - & \frac{\tilde{a}_1}{x^{2\alpha+1}} - \frac{a_2}{x^{2\alpha+1}}=o_{+\infty} \left(x^{-\beta}\right). \\
        & \text{Case } p>\frac{2\alpha+1}{\alpha+1}: \quad \exists \beta> 2\alpha+1, & Q(x)- \frac{a_1}{x^{\alpha+1}} - &   \frac{a_2}{x^{2\alpha+1}} =o_{+\infty} \left(x^{-\beta}\right).
    \end{aligned}
\end{equation*}
\end{theo}

\begin{rema}
Let us notice that the constants involved in the asymptotic expansion are coherent with other situations. In the case $\alpha=1$ and $f(u)=u^2$, thus for (BO), only the terms with an even power are present, since $a_2=0$, but $a_1\neq 0$. On the other hand, the case $\alpha=2$ and $f(u)=u^p$ with $p$ an integer corresponds to the generalized Korteweg-de Vries equation whose solitons have an exponential decay. By replacing $\alpha=2$ in the coefficients of the asymptotic expansion, we find $a_1=a_2=0$, which is coherent with the exponential decay.
\end{rema}

The next theorems refine the asymptotic development of $Q$ in the case of a polynomial non-linearity.

\begin{theo}\label{propo:poly_derivatives}
Let $p\in\N$, $p\geq2$, $\alpha\in\left(\frac{p-1}{1+p},2\right)$, and $Q$ be solution of $|D|^{\alpha}Q+Q-Q^{p}=0$ verifying condition \eqref{hypo:Q}. Then $Q\in H^{\infty}(\R)$ and verifies that:
    \begin{align*}
        \forall j\in\N, \quad \exists \beta=\beta(j)> \alpha+1 +j, \quad Q^{(j)}(x) -(-1)^{j}\frac{(\alpha+j)!}{\alpha!} \frac{k_1}{x^{\alpha+1+j}}\int Q^p  =o_{+\infty} \left(x^{-\beta}\right).
    \end{align*}
\end{theo}

In the case of a cubic non-linearity, the next theorem provides a sharper asymptotic development.

\begin{theo}[Higher order expansion]\label{thm:asympQ}
Let $\alpha\in(1,2)$, $p=3$, and $Q$ be a solution of $|D|^{\alpha}Q+Q-Q^{p}=0$ verifying condition \eqref{hypo:Q}. Then, there exists a constant $C=C(\alpha,p)>0$ such that:
\begin{align}
    & \left\vert Q(x) -\left( \frac{a_1}{x^{\alpha+1}} + \frac{a_2}{x^{2\alpha+1}} + \frac{a_3}{x^{\alpha+3}} \right)\right\vert \leq \frac{C}{x^{3\alpha+1}} ,\label{asympt:Q}\\
    & \left\vert Q'(x) +(\alpha+1)\frac{a_1}{x^{\alpha+2}} + (2\alpha+1) \frac{a_2}{x^{2\alpha+2}} \right\vert  \leq \frac{C}{x^{3\alpha+1}}. \label{asympt:Q'}
\end{align}
\end{theo}

\begin{rema}
The constants obtained in the asymptotic development of $Q(x)$ are dependent of the functions $x^lk(x)$ and $x^l|Q|^{p-1}(x)Q(x)$, for $l\in\N$, which are dependent on $p$ and $\alpha$ as for $a_3$ defined in \eqref{defi:a1_a2_a3}. Thus, our method allows a further asymptotic development while $x^l|Q|^{p-1}(x)Q(x)$ and $x^lk(x)$ remain integrable.   
\end{rema}

Let us point out an application of these results. The authors of this article describe in \cite{EV22} the long term interaction of two solitary waves with the same velocity for the \ref{fgKdV} equation. To this aim, they used the asymptotic behaviour of those waves to quantify the distance between the two objects. Indeed, the asymptotic behaviour up to order $3$ given in Theorem \ref{propo:poly_derivatives} was necessary to quantify the strong interaction.

A natural question is the generalization of those results to higher dimensions. The authors do not know how the computations introduced in this article can be generalized for higher dimensions, and if similar theorems can be stated in this new setting.

\subsection{Ideas of the proof}
Let us describe the main ideas to obtain the asymptotic developments. We use the kernel formulation, introduced first by Bona-Li \cite{BL97}, of the equation satisfied by $Q$: $Q=k\star f(Q)$. Let us explain formally why the asymptotic behaviour of $Q$ is characterised by the one of $k$, as remarked by Bona-Li \cite{BL97}. Consider the general formulation:
\begin{align}
    -cu + \mathcal{L}u + f(u)=0, \quad u:\R^n \rightarrow \R,
\end{align}
where $\mathcal{L}$ stands for a linear operator with symbol $m$, and $f$ a general non-linearity satisfying $\vert f (u)\vert \simeq \vert u \vert ^p$ for $p>1$. Finding a solution to this equation is equivalent to find $u$ of the form:
\begin{align*}
    u(x) =\int_{\R^n} k(y-x) f(u(y)) dy,
\end{align*}
where the kernel $k$ is given by:
\begin{align}\label{defi:k_fourier}
    k(x):=\mathcal{F}^{-1}\left(\frac{1}{c-m(\xi)}\right)(x),
\end{align}
where $\mathcal{F}^{-1}$ stands for the inverse Fourier transform. Let us suppose that $Q$ decays at infinity. The asymptotic behaviour of $Q$ at infinity is given by the largest term in the asymptotic behaviours of $k$ and of $f(Q)$ at infinity. Indeed, for $x$ large, from the convolution formula, $k(x-\cdot)$ and $f(Q)$ are localized at different places.
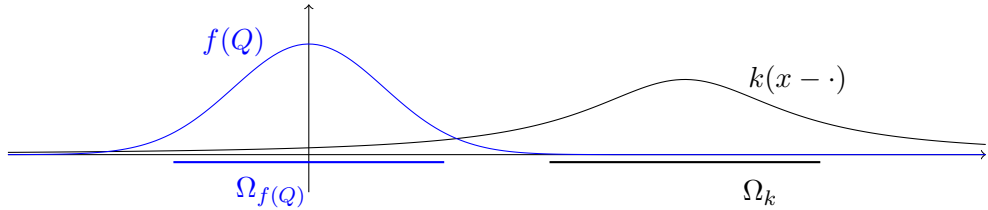
\begin{center}
\begin{tikzpicture}
    \draw [->] (0,-0.5)-- (0,2);
    \draw [->] (-4,0)-- (9,0);

    \draw[blue] (-1,1.5) node {$f(Q)$};
    \draw [domain=-4:9, samples=200] plot(\x,{1/(1+(\x-5)*(\x-5)/2.5)});
    \draw (6.5,1) node {$k(x-\cdot)$};
    \draw [domain=-4:9, samples=200, color=blue] plot(\x,{4*exp(-(1+\x*\x/2))});

    \draw[blue] (-0.5,-0.5) node {$\Omega_{f(Q)}$};
    \draw[thick, blue] (-1.8,-0.1)--(1.8,-0.1);
    \draw (6,-0.5) node {$\Omega_k$};
    \draw[thick] (3.2,-0.1)--(6.8,-0.1);
\end{tikzpicture}
\end{center}
One can notice that on $\Omega_{f(Q)}$, where the mass of $f$ is located, the tail of $k(x-\cdot)$ is larger than the tail of $f(Q)$ on $\Omega_k$, where the mass of $k(x-\cdot)$ is located. By the definition of $\vert f (u)\vert \simeq \vert u \vert ^p$, with $p>1$, the decay of $f(Q)$ is larger than the one $Q$. In other words, the main order term in the integral formulation is given by the decay of $k$ on $\Omega_{f(Q)}$.

The first part of the article is dedicated to the study of integrability and of a potential differentiability of $k$. The key point of the analysis of $k$ is the use of complex analysis techniques developed in the beginning of the twentieth century. Following the road map of P\'olya \cite{Pol23} developed in 1923 (see also Blumenthal-Getoor \cite{BG60}), we rewrite the function $k$ with the help of an auxiliary function $h$. The asymptotic development of $h$ is obtained by complex analysis tools : we consider the holomorphic extension of $h$ in different regions to extract the terms at the main orders. We conclude that $k$ is in $L^1(\R)$ for any $\alpha\in(0,2)$, and for $\alpha\in (1,2)$ the derivative of $k$ exists and is also in $L^1(\R)$. The bottleneck of the study of $k$ stands in its behaviour at $0$: having a $L^1$-bound of this operator at $0$ is sufficient to get the asymptotic behaviour of $Q$.

Once the asymptotic development of $k$ is established, we use this development to obtain the one of $Q$ with the previous convolution formula. Since the assumption of $Q$ only consists in a slow decay, we first improve the property of decay by injecting successively this bound into the convolution. The limit point of improvement of the decay is when the biggest term of the asymptotic expansion is coming from the asymptotic development of $k$. We obtain the main order term in the development of $Q$ as being the same as the one of $k$. Then, to get the asymptotic development of the first derivative of $Q$, we need the differentiability of $k$, and our method applies only if $\alpha>1$. By using the same arguments as before, we obtain the asymptotic at first order of the derivatives of $Q$. Similarly, we obtain the first order development of the derivatives of $Q$, while those derivatives exist. The method exposed in this paragraph has been exploited by Chen-Walsh-Wheeler \cite{CWW19} to derive the asymptotic expansion of deep water solitary waves with localized vorticity in dimension 2 and 3.

The previous results can be extended if $Q$ is a positive solution. Indeed, the previous issue of the existence of a finite number of derivatives was due to the non-regularity of $u \mapsto \vert u\vert^{p-1} u$ applied at $Q$, so potentially applied at $0$. Nevertheless the assumption of $Q$ positive circumvents this issue by studying the regularity of the function $u \mapsto \vert u\vert^{p-1} u$ on $\R_+^*$ only. This function is smooth, and so is $Q$. Then, we notice that having $Q$ positive also implies that $Q$ is an even function, and this property helps to get the next term in the asymptotic of $Q$. We distinguish the different cases of balance between the non-linearity and the dispersion. If the non-linearity is to large, the next order term comes from the development of $k$; if the non-linearity is to small, the next order term comes from the development of $Q^p$; in the case of exact balance, the two previous terms are at the same order and the sum of them furnishes the next order term.

Finally, in Theorem \ref{thm:asympQ} we study the equation in the case of a particular non-linearity, and the asymptotic development at order $3$, is obtained  originating only from the asymptotic development of $k$.

\subsection{Related results}

For many equations, the solitary waves can be defined by a convolution formula, and their asymptotic expansion is a consequence of the regularity of the kernel. From the definition of $k$ by the Fourier transform, see \eqref{defi:k_fourier}, the decay of $k$ is given by the regularity of $(c-m)^{-1}$. Thus, if $(c-m)^{-1}$ is a smooth function, we expect to get $Q$ exponentially decaying, whereas if $(c-m)^{-1}$ is not smooth the solution $Q$ should vanish at infinity algebraically. Let us give an overview of those phenomena for diverse equations.

We begin with exponentially decaying solutions. For the local dispersion setting, the generalized KdV equation admits solitary waves which decay exponentially, with the kernel $k=\mathcal{F}^{-1}\left( \frac{1}{c+|\xi|^2} \right)$. Some non-local equations admit solitary waves with exponential decay, like the Whitham type equations as studied by Bruell-Ehrnstr{\"o}m-Pei and Arnesen \cite{BEP17,Arn22}. For solitary waves with velocity $c>1$, the symbol of dispersion of the Whitham equation is given by $m(\xi)=\left( \frac{\tanh(\xi)}{\xi} \right)^{\frac{1}{2}}$, and thus $(c-m(\xi))^{-1}$ is not in $L^2(\R)$. Even though $k$ is not well-defined, by an elegant reformulation in \cite{BEP17,Arn22}, the study of the asymptotic behaviour of the solitary waves is equivalent to the study of the kernel operator $\tilde{k}(x):= \mathcal{F}^{-1}(\frac{m}{c-m})(x)$. Because this last formulation is given by a smooth operator $\frac{m}{c-m}$ in $L^2(\R)$, \cite{BEP17} proved the exponential decay of the solitary waves. Following the lines of \cite{BEP17}, Pei in \cite{Pei20} proved the exponential decay of the solitary waves for the Degasperis-Procesi equation.

Other equations are known to own solitary waves with algebraic decay. In the local setting, the nonlinear wave equation in dimension $3\leq d\leq 5$, with the kernel $k=\mathcal{F}^{-1}\left(\vert \xi \vert^{-2}\right)$ (defined in the weak sense), admits steady waves with algebraic decay proved by Gidas-Ni-Nirenberg \cite{GNN79}. For the non-local case one can cite the generalized Benjamin-Ono equation studied by Mari\c{s} \cite{Mar02}, with the kernel $k=\mathcal{F}^{-1}\left(\frac{1}{c+|\xi|} \right)$, or the fKdV equation investigated  by Frank-Lenzmann-Silvestre \cite{FLS16}, where the kernel is given by $k=\mathcal{F}^{-1}\left(\frac{1}{c+|\xi|^{\alpha}} \right)$, with $0<\alpha<2$.

For non-radial solutions in dimensions larger than $2$, the asymptotic behaviour when $\vert x \vert$ tend to $+\infty$ can depend on the direction. Indeed, in de Bouard-Saut \cite{dBS97}, an algebraic bound of the asymptotic behaviour is obtained for a solitary wave of the generalized Kadomstev-Petviashvili equation; furthermore, there exists one direction for which this bound is optimal. This kind of asymptotic behaviour has also been observed by Gravejat \cite{Gra06}: by studying the solitary waves of the Gross-Pitaevskii equation, he gave an algebraic asymptotic behaviour of order one, whose coefficient at the main order depends explicitly on the angle. For the Benjamin-Ono--Zakharov-Kuznetsov equation, Esfahani-Pastor-Bona \cite{BEP15} proved that the solitary waves of this $2$-dimensional equation decay at least polynomially in one direction, and faster than any polynomial in another direction. The question of asymptotic behaviour for several dimensions is thus more difficult to answer, in particular if the dispersion operator is not rotationally invariant.

When the non-linearity is not an algebraic function, the asymptotic behaviour can depend at the main order on the dispersion operator and on the non-linearity. For the non-local Gross-Pitaevskii equation, de Laire-L{\'o}pez-Mart{\'i}nez \cite{dLL22} described the asymptotic behaviour of solitary waves $v_c$ depending on a convolution function in the non-linearity. The asymptotic decay of $1-\vert v_c \vert^2$ can be algebraic or exponential at infinity depending of the choice of the convolution function.

The KdV type equations are toy models for some regimes of water waves. The Euler equation with a free surface give a more accurate description of those waves. Following ideas introduced by Li-Nirenberg \cite{LN98}, Ehrnstr{\"o}m-Walsh-Zeng \cite{EWZ22} studied stationary capillary-gravity waves in a two-dimensional body of water that rests above a flat ocean bed and below vacuum, which is related to the existence of solutions of semi-linear equations. They proved the existence of exponentially decaying stationary waves with a non-trivial vorticity. In the infinitely deep body of water, Chen-Walsh-Wheeler \cite{CWW19} investigated the existence of capillary-gravity solitary waves with compactly supported vorticity: they proved their existence under some additional assumptions on the vorticity, and proved the algebraic decay of those waves.

\subsection{Outline of the paper}

The second section is dedicated to the study of the kernel $k$. Section 3 is dedicated to the proof, in the general case with no assumption on $Q$, $\alpha$ nor $p$, of the the first order asymptotic of $Q$, see Theorem \ref{propo:notre_method_1}. Section 4 deals with the proof of Theorem \ref{propo:notre_method_2} with the case of more regularity. Finally, one particular case is dealt with in Section 5 where the asymptotic expansion is given at order $3$, see Theorem \ref{thm:asympQ}. The Appendix recalls the proof of the regularity of $Q$ if $p$ is an integer.

\subsection{Notations}

The japanese bracket $\langle \cdot \rangle$ is defined on $\R$ by $\langle x \rangle := (1+\vert x \vert^2)^\frac{1}{2}$.

If $\Omega$ is a subset of $\R$, we denote by $\mathcal{C}^k(\Omega)$ the set of $k$-differentiable functions, with the usual generalization for $k=\infty$.

By denoting $\lambda$ the Lebesgue measure, we define the generalized $L^p$ space as:
\begin{align*}
    \forall p \in (1,+\infty), \quad L^p(\mathbb{R}) := \left\{ f \text{ function on } \Omega, \quad \Omega \subset \mathbb{R} \text{ and } \lambda( \Omega^C)=0, \quad \| f\|_{L^p(\Omega)}<\infty \right\}.
\end{align*}
By abuse of notations, we denote by $\|\cdot\|_{L^p}$ the $L^p$-norm over any subset $\Omega$ of $\R$ over which the function is well-defined, with $\lambda(\Omega^C)=0$.

We denote the Fourier transform by $\mathcal{F}$, defined by:
\begin{align*}
    \forall u \in L^2(\R), \quad \mathcal{F}(u)(\xi) := \int_\R e^{-2i\pi x\xi} u(x) dx,
\end{align*}
and its inverse by $\mathcal{F}^{-1}$. The Sobolev space for $s\in \R$ defined by the Fourier transform is thus:
\begin{align*}
    H^s(\R) := \left\{ u \in L^2(\R); \int_\R \langle \xi \rangle^{2s} \left\vert \mathcal{F}(u)(\xi) \right\vert^2 d\xi<\infty \right\}.
\end{align*}

The usual convolution operator $\star$ is defined by:
\begin{align*}
    \forall f,g \in L^2(\R), \quad f\star g(x) := \int_\R f(x-y) g(y)dy.
\end{align*}

For a fixed $x\in\R$, we define the following subset of $\mathbb{R}$:
\begin{align}\label{defi:omega_x}
    \Omega_x:= \left\{ \vert y \vert \leq \frac{x}{2}  \right\}.
\end{align}

We denote by $C$ a positive constant that can change for line to line.

\section{Asymptotic expansion of the kernel operator}

Let us define the kernel associated to the operator $(1+\vert D \vert^\alpha)^{-1}$:
\begin{align}\label{defi:k}
k(x):=\mathcal{F}^{-1}\left(\frac{1}{1+|\xi|^{\alpha}}\right)(x).  
\end{align}

For $\alpha \in (0,1]$, the inverse Fourier transform of $(1+\vert \xi \vert^\alpha)^{-1}$ is understood as an improper integral for each $x\in \mathbb{R}^*$.

In order to give the asymptotic of the ground state, we establish an asymptotic development of $k$. Lemma \ref{lemma:est:invlapl} is dedicated to this asymptotic development. As an application, we prove in Corollary \ref{lemma:est:kernel} that some decay properties are preserved under the convolution with $k$. 

\begin{lemm}\label{lemma:est:invlapl} The function $k$ is in $\mathcal{C}^\infty(\R^{*})$, and in $L^1(\R)$. Furthermore, there exists a sequence $(k_n)_n \in \mathbb{R}^\mathbb{N}$, such that for any $N>0$, there exists $C_N>0$ such that:
\begin{align}\label{est:invlapl}
   \forall \vert x \vert >1, \quad \left|k(x)- \sum_{n=1}^N \frac{k_n}{ \vert x\vert ^{n\alpha+1}}  \right|\leq \frac{C_N}{\vert x\vert ^{(N+1)\alpha+1}},
\end{align}
and $k'$ admits the following development:
\begin{align}\label{est:invlapl_deriv}
    \forall \vert x\vert  >1, \quad \left\vert k'(x) - \sign(x) \sum_{n=1}^N \frac{(n\alpha+1)k_n}{ \vert x \vert^{n\alpha+2}} \right\vert \leq \frac{C_N}{ \vert x\vert  ^{(N+1)\alpha+2}}.
\end{align}
Furthermore, in the case $\alpha>1$ $k'\in L^1(\R)$.
\end{lemm}

Notice that the values of $k_1$ and $k_2$ are defined in \eqref{defi:k1_k2}.

\begin{rema}
Frank, Lenzmann and Silvestre established in \cite{FLS16}, Lemma C.1., the asymptotic of the kernel at the first order, and proved that the kernel $k$ is positive.
\end{rema}

\begin{rema}
The previous development only holds for fixed $N$, we ignore if the serie converges.     
\end{rema}

Let us first write this kernel with a more convenient formula, by defining the function $h$ by:
\begin{align}
    h(y)= \int_0^\infty \cos\left( y \eta \right) e^{-\eta^\alpha} d\eta.
\end{align} 
We claim the formula:
\begin{claim}\label{claim:k_function_of_h}
The following equality holds:
    \begin{align*}
      \forall x \neq 0, \quad k(x) = \frac{1}{\pi} \int_0^{\infty} \frac{e^{-s}}{s^\frac{1}{\alpha}} h\left( \frac{x}{s^\frac{1}{\alpha}} \right) ds.
    \end{align*}
\end{claim}

One can notice that this formula is more convenient. Indeed, for any value of $\alpha\in(0,2)$, the asymptotic development of $h$ at infinity induces that $k$ is $\mathcal{C}^\infty$ on $\mathbb{R}^*$, whereas this property was not clear by the definition by the Fourier transform. Notice that Lemma C.1 of \cite{FLS16} used another formulation of $k$ : by applying Bernstein's theorem on the fractional heat kernel, $k$ is rewritten by the semi-group associated with the heat kernel and an adequate non-negative finite measure.

Let us begin by proving Claim \ref{claim:k_function_of_h}. To do so, we need the asymptotic development of $h$:

\begin{claim}\label{claim:h_asymptotic_dvt}
The following expansions hold:
\begin{align}\label{estimation:h}
    \forall y>1, \quad \left|h(y) -\frac{\pi k_1}{y^{\alpha+1}}- \frac{\pi k_2}{2y^{2\alpha+1}} \right| \leq \frac{C}{y^{3\alpha+1}}, 
\end{align}
and
\begin{align}\label{estimation:h_prime}
    \forall y>1, \quad \left|h'(y) -  \frac{k_1'}{y^{\alpha+2}} \right| \leq \frac{C}{y^{2\alpha+2}}, \quad k_1'\in\R.
\end{align}
Furthermore, for any $\beta \in (0,1]$, there is another constant $C=C(\beta)$ such that:
\begin{align}\label{eq:bound_h}
    \forall \vert y \vert \leq 1, y \neq 0,  \quad \vert yh(y) \vert \leq Cy^{\beta} \quad \text{and} \quad \vert h'(y) \vert \leq \frac{C}{y}.
\end{align}
\end{claim}

\begin{proof}[Proof of Claim \ref{claim:k_function_of_h}]
In the case $\alpha\in (1,2)$, the proof holds by a Fubini's argument, as in \cite{BG60}. We give the general proof that also holds in the case $\alpha\in(0,1]$.

Using that $\displaystyle\frac{1}{t}= \int_{\mathbb{R}^+}e^{-st}ds$ :
\begin{align*}
k(x) = \frac{1}{2\pi} \int_{\mathbb{R}}e^{ix\xi} \int_0^1 e^{-s \vert \xi \vert^\alpha} e^{-s} ds d\xi +\frac{1}{2\pi} \int_{\mathbb{R}}e^{ix\xi} \int_1^\infty e^{-s \vert \xi \vert^\alpha} e^{-s} ds d\xi.
\end{align*}

We then apply Fubini's theorem for the first integral. By using Claim \ref{claim:h_asymptotic_dvt}, it is also possible to use Fubini's theorem on the first integral: the asymptotic development of $h$ is necessary, in particular if $\alpha\in(0,1]$. Finally, the change of variable $\eta^\alpha = \xi^\alpha s$ concludes Claim \ref{claim:k_function_of_h}.
\end{proof}

\begin{proof}[Proof of Claim \ref{claim:h_asymptotic_dvt}]
The next ideas are inspired from P\'olya \cite{Pol23}. By integration by parts and changing the variable $u^{\frac{1}{\alpha}}=y\eta$, we deduce that:
\begin{align*}
    y^{1+\alpha}h(y)=y^{\alpha}\int_{0}^{+\infty}\frac{d}{d\eta}(\sin(y\eta))e^{-\eta^{\alpha}}d\eta =\int_{0}^{+\infty}\sin\left(u^{\frac{1}{\alpha}}\right)e^{-\frac{u}{y^{\alpha}}}du=\text{Im}\left(\int_{0}^{+\infty} e^{iu^{\frac{1}{\alpha}}-\frac{u}{y^{\alpha}}}du\right).
\end{align*}
Note that the previous integral is not well defined for $y=+\infty$. To bypass this difficulty, we apply contour integration, see \cite{Pol23}. Let: 
\begin{align*}
D_n:=\{re^{i\frac{\pi}{4}\alpha}: r\in(0,n]\}, \quad D:=\{re^{i\frac{\pi}{4}\alpha}: r\in(0,+\infty)\},\quad C_n:=\{ne^{i\gamma}:\gamma\in\left[0,\frac{\pi}{4}\alpha\right] \}.
\end{align*} 
We set $\gamma_n$ the curves that range $D_n$, $C_n$ and then $[0,n]$ counterclockwise. Since $u \mapsto e^{iu^{\frac{1}{\alpha}}-\frac{u}{y^{\alpha}}}$ is holomorphic in $\R^*_{+}+i\R_{+}$, we deduce that:
\begin{align*}
    0=\lim_{n\to+\infty}\int_{\gamma_n}e^{iu^{\frac{1}{\alpha}}-\frac{u}{y^{\alpha}}}du=\int_{D}e^{iu^{\frac{1}{\alpha}}-\frac{u}{y^{\alpha}}}du - \int_{0}^{+\infty}e^{iu^{\frac{1}{\alpha}}-\frac{u}{y^{\alpha}}} du.
\end{align*}
Furthermore, we get that:
\begin{align*}
    y^{1+\alpha}h(y)= \text{Im}\left( \int_{0}^{+\infty}\exp\left(\frac{\sqrt{2}}{2}(i-1) r^{\frac{1}{\alpha}} - \frac{re^{i\frac{\pi}{4}\alpha}}{y^{\alpha}}   \right) e^{i\frac{\pi}{4}\alpha} dr \right).
\end{align*}
To obtain the asymptotic expansion of $y^{1+\alpha}h(y)$ at $+\infty$, we split the former integral in two parts. We set: 
\begin{align*}
    J_1:=\text{Im}\left(\int_{0}^{y^{\frac{\alpha}{2}}}\exp\left(\frac{\sqrt{2}}{2}(i-1) r^{\frac{1}{\alpha}} - \frac{re^{i\frac{\pi}{4}\alpha}}{y^{\alpha}}   \right) e^{i\frac{\pi}{4}\alpha} dr \right),\\ J_2:=\text{Im}\left(\int_{y^{\frac{\alpha}{2}}}^{+\infty}\exp\left(\frac{\sqrt{2}}{2}(i-1) r^{\frac{1}{\alpha}} - \frac{re^{i\frac{\pi}{4}\alpha}}{y^{\alpha}}   \right) e^{i\frac{\pi}{4}\alpha} dr\right).
\end{align*}
Since $\alpha\in(1,2)$, we have that:
\begin{align}\label{estimation:J2}
    |J_2|\leq \int_{y^{\frac{\alpha}{2}}}^{+\infty} \exp\left( -\frac{\sqrt{2}}{2}r^{\frac{1}{\alpha}} - \cos\left(\frac{\pi}{4}\alpha \right)\frac{r}{y^{\alpha}}\right)dr \leq C \int_{y^{\frac{\alpha}{2}}}^{+\infty} \exp\left( -\frac{\sqrt{2}}{2}r^{\frac{1}{\alpha}} \right) dr \leq Ce^{-\frac{\sqrt{y}}{4}}.
\end{align}
Now, we estimate $J_1$. First, we rewrite $J_1$ as: 
\begin{align}
    J_1=&\text{Im}\left(\int_{0}^{y^{\frac{\alpha}{2}}}\exp\left(i \left(re^{i\frac{\pi}{4}\alpha}\right)^{\frac{1}{\alpha}} - \frac{re^{i\frac{\pi}{4}\alpha}}{y^{\alpha}}   \right) e^{i\frac{\pi}{4}\alpha} dr \right) \notag\\
    =&\text{Im}\left(\int_{0}^{y^{\frac{\alpha}{2}}}\exp\left(i \left(re^{i\frac{\pi}{4}\alpha}\right)^{\frac{1}{\alpha}} +i\frac{\pi}{4}\alpha\right) \left[  \exp\left(-\frac{re^{i\frac{\pi}{4}\alpha}}{y^{\alpha}}\right) - 1 + \frac{re^{i\frac{\pi}{4}\alpha}}{y^{\alpha}} 
    \right]  dr \right) \notag\\
    +& \text{Im}\left(\int_{0}^{y^{\frac{\alpha}{2}}}\exp\left(i \left(re^{i\frac{\pi}{4}\alpha}\right)^{\frac{1}{\alpha}} +i\frac{\pi}{4}\alpha\right) \left[ 1 - \frac{re^{i\frac{\pi}{4}\alpha}}{y^{\alpha}} 
    \right]  dr \right)=:J_{11} + J_{12}. \label{eq:decomposition_J_1}
\end{align}
From the Taylor expansion of $e^z$, we deduce that:
\begin{align}\label{estimation:J11}
    |J_{11}|\leq C\int_{0}^{y^{\frac{\alpha}{2}}}e^{-\frac{\sqrt{2}}{2}r^{\frac{1}{\alpha}}}\left(\frac{r}{y^{\alpha}}\right)^2 dr\leq \frac{C}{y^{2\alpha}}.
\end{align}
Let us rewrite $J_{12}$:
\begin{align}\label{reecrire:J12}
    J_{12}&=\text{Im}\left(\int_{0}^{+\infty}e^{\frac{\sqrt{2}}{2}(i-1) r^{\frac{1}{\alpha}}+i\frac{\pi}{4}\alpha} dr\right) -\text{Im}\left( \int_{y^{\frac{\alpha}{2}}}^{+\infty}e^{\frac{\sqrt{2}}{2}(i-1) r^{\frac{1}{\alpha}}+i\frac{\pi}{4}\alpha} dr\right)  \\
    &-\frac{1}{y^{\alpha}}\text{Im}\left(\int_{0}^{+\infty}re^{\frac{\sqrt{2}}{2}(i-1) r^{\frac{1}{\alpha}}+i\frac{\pi}{2}\alpha} dr\right) +\text{Im}\left( \frac{1}{y^{\alpha}}\int_{y^{\frac{\alpha}{2}}}^{+\infty}re^{\frac{\sqrt{2}}{2}(i-1) r^{\frac{1}{\alpha}}+i\frac{\pi}{2}\alpha} dr\right)\notag.
\end{align}
Arguing similarly as \eqref{estimation:J2}, we get that:
\begin{align}\label{estimation:J12queue}
    \left|\int_{y^{\frac{\alpha}{2}}}^{+\infty}e^{\frac{\sqrt{2}}{2}(i-1) r^{\frac{1}{\alpha}}+i\frac{\pi}{4}\alpha} dr \right| + \left| \frac{1}{y^{\alpha}}\int_{y^{\frac{\alpha}{2}}}^{+\infty}re^{\frac{\sqrt{2}}{2}(i-1) r^{\frac{1}{\alpha}}+i\frac{\pi}{2}\alpha} dr \right| 
    \leq C e^{-\frac{\sqrt{y}}{4}}. 
\end{align}
In order to get the sign of the two integrals over $\mathbb{R}^{+}$ in \eqref{reecrire:J12}, we use again a contour integration. We set:
\begin{align*}
\widetilde{D_n}:=\{re^{i\frac{\pi}{2}\alpha}: r\in(0,n]\}, \quad \widetilde{D}:=\{re^{i\frac{\pi}{2}\alpha}: r\in(0,+\infty)\},\quad \widetilde{C_n}:=\{ne^{i\gamma}:\gamma\in\left[\frac{\pi}{4}\alpha,\frac{\pi}{2}\alpha\right] \}.
\end{align*} 
Therefore, by using the curves who range $\widetilde{D_n}$, $\widetilde{C_n}$ and $D_n$ counterclockwise, we have that:
\begin{align}\label{contour2}
    &\text{Im}\left(\int_{0}^{+\infty}e^{\frac{\sqrt{2}}{2}(i-1) r^{\frac{1}{\alpha}}+i\frac{\pi}{4}\alpha} dr\right)=\text{Im}\left(\int_{0}^{+\infty}e^{- r^{\frac{1}{\alpha}}+i\frac{\pi}{2}\alpha} dr\right) = \sin(\frac{\pi}{2}\alpha) \Gamma(\alpha+1) ,\\
    &\text{Im}\left(\int_{0}^{+\infty}re^{\frac{\sqrt{2}}{2}(i-1) r^{\frac{1}{\alpha}}+i\frac{\pi}{2}\alpha} dr\right)=\text{Im}\left(\int_{0}^{+\infty}re^{- r^{\frac{1}{\alpha}}+i\pi\alpha} dr\right) = \frac{1}{2} \sin(\pi \alpha) \Gamma(2\alpha+1).  \notag
\end{align}
Gathering \eqref{estimation:J2}, \eqref{estimation:J11}, \eqref{estimation:J12queue} and \eqref{contour2} we conclude \eqref{estimation:h}, with $k_1$ and $k_2$ recalled in \eqref{defi:k1_k2}.

We continue by proving \eqref{estimation:h_prime}. 
The proof of this property is also obtained by contour argument as for $h$, Claim \ref{claim:h_asymptotic_dvt}. Let us give some steps of the proof. By integration by part, and the change of variables $u=(y\eta)^2$ we obtain that:
\begin{align*}
    y^{2}h'(y)=-\frac{1}{2}\int_{0}^{\infty}\sin\left(u^{\frac{1}{2}} \right) e^{-u^{\frac{\alpha}{2}}y^{-\alpha}}du.
\end{align*}
First, using the contour defined by  
\begin{align*}
D_n':=\left\{re^{i\frac{\pi}{2}}: r\in(0,n]\right\}, \quad D':=\left\{re^{i\frac{\pi}{2}}: r\in(0,+\infty)\right\},\quad C_n':=\left\{ne^{i\gamma}:\gamma\in\left[0,\frac{\pi}{2}\right] \right\},
\end{align*} 
we have
\begin{align*}
    y^2h'(y)= -\frac{1}{2} \text{Im} \left( \int_0^{+\infty} \exp \left( \sqrt{r} e^{i\frac{3}{4}\pi} - \frac{r^\frac{\alpha}{2}}{y^\alpha} e^{i\frac{\pi}{4}\alpha} +i\frac{\pi}{2} \right) dr \right).
\end{align*}
By using a decomposition as in \eqref{estimation:J2} and \eqref{eq:decomposition_J_1} with the bound of the different terms, we have:
\begin{align*}
    \left\vert y^2 h'(y) + \frac{1}{2} \text{Im} \left( \int_0^{+\infty} \exp  \left( (re^{i\frac{3}{2}\pi})^{\frac{1}{2}} +i\frac{\pi}{2} \right) \left( 1 - \frac{r^{\frac{\alpha}{2}}}{y^\alpha}e^{i\frac{\pi}{4}\alpha} \right) dr \right)\right\vert \leq \frac{C}{y^{2\alpha}}.
\end{align*}
Using a second contour integration defined by 
\begin{align*}
\widetilde{D_n'}:=\left\{re^{i\frac{3\pi}{2}}: r\in(0,n]\right\}, \quad \widetilde{D}':=\left\{r: r\in(0,+\infty)\right\},\quad \widetilde{C_n'}:=\left\{ne^{i\gamma}:\gamma\in\left[\frac{3\pi}{2},2\pi\right] \right\},
\end{align*} 
and the fact 
\begin{align*}
    \Re\left(\int_{0}^{+\infty} e^{\frac{\sqrt{2}}{2} (i-1)\sqrt{u}} du\right)=2\Re\left(\int_{0}^{+\infty} ve^{\frac{\sqrt{2}}{2} (i-1)v} dv\right)=0,
\end{align*}
we obtain the asymptotic development \eqref{estimation:h_prime} of $h'$.

To prove \eqref{eq:bound_h}, we multiply $h$ by $y^{1-\beta}$ and integrate by part:
\begin{align*}
    y^{1-\beta} h(y)= \int_0^{+\infty} \frac{\sin(y\eta)}{(y\eta)^\beta} \alpha \eta^{\alpha-1+\beta} e^{-\eta^\alpha} d\eta.
\end{align*}
Since $\frac{\sin(z)}{z^\beta}$ is uniformly bounded in $\beta \leq 1$, we obtain the first inequality of \eqref{eq:bound_h}. For the second inequality:
\begin{align*}
    y^2h'(y)= - y -y \int_0^{+\infty} \frac{\sin(y\eta)}{y\eta} \alpha \left( (\alpha-1) \eta^{\alpha-1} - \alpha \eta^{2\alpha-1} \right) e^{-\eta^\alpha} d\eta,
\end{align*}
and a direct bound gives the second part of \eqref{eq:bound_h}.
\end{proof}

We continue with the proof of Lemma \eqref{lemma:est:invlapl}. We focus on the case $N=2$, the general statement is obtained using the same proof.

\begin{proof}[Proof of Lemma \ref{lemma:est:invlapl}]
Concerning the regularity of $k$, our description of the kernel is not well-suited to get this regularity. However, in Lemma C.1 of \cite{FLS16}, the kernel $k$ is proved to be in $\mathcal{C}^{\infty}(\mathbb{R}^*)$.

We continue with the asymptotic expansion of $k$ at infinity. By parity of $k$, we focus on the case $x>0$. By Claim \ref{claim:k_function_of_h}, we decompose the integral according to the values of $s$; in other words, we need to study carefully when $\frac{x}{s^\frac{1}{\alpha}}$ tends to $+\infty$ as $x\rightarrow + \infty$. We thus decompose:
\begin{align*}
   \pi k (x) =\int_0^{x^{\frac{\alpha}{2}}} \frac{e^{-s}}{s^{\frac{1}{\alpha}}} h\left( \frac{x}{s^{\frac{1}{\alpha}}} \right) ds +\int_{x^{\frac{\alpha}{2}}}^{+\infty} \frac{e^{-s}}{s^{\frac{1}{\alpha}}} h\left( \frac{x}{s^{\frac{1}{\alpha}}} \right) ds =I_1+I_2.
\end{align*}

Let us first find a bound of $I_2$. Since $h$ is a bounded function, we get:
\begin{align*}
    \vert I_2 \vert \leq Ce^{-x^\frac{\alpha}{4}}.
\end{align*}

The asymptotic expansion of $I_1$ is now given by the the one of $h$, since $\frac{x}{s^\frac{1}{\alpha}} \geq x^\frac{1}{2}$ for $s \leq x^{\frac{\alpha}{2}}$.
Since $\displaystyle\int_{\mathbb{R}}se^{-s}ds=1$, $\displaystyle\int_{\mathbb{R}}s^2e^{-s}ds=2$ and by using \eqref{estimation:h}, we conclude the proof of the asymptotic development of k. 

This development justifies the integrability of $k$ at $\infty$. We now justify the integrability of $k$ around $0$:
\begin{align*}
    xk(x) = \frac{1}{\pi} \int_0^{x^\alpha} e^{-s} \left( yh(y) \right)_{\vert y =xs^{-\frac{1}{\alpha}}} ds + \frac{1}{\pi} \int_{x^\alpha}^{+\infty} e^{-s} \left( yh(y) \right)_{\vert y =xs^{-\frac{1}{\alpha}}} ds.
\end{align*}
For the first integral, we use \eqref{estimation:h} for $\vert y \vert \geq 1$, and the first integral is bounded by $x^\alpha$. For the second integral, we use \eqref{eq:bound_h} for $\vert y \vert \leq 1$ with $\beta=\min(1,\alpha)$ and we compute the integral:
\begin{align*}
    \text{if } \alpha \leq 1, \quad \int_{x^\alpha}^{+\infty} e^{-s} \left( \frac{x}{s^\frac{1}{\alpha}} \right)^{\alpha} ds \leq C x^{\alpha} \vert \ln(x) \vert;\quad \text{if } \alpha >1, \quad \int_{x^\alpha}^{+\infty} e^{-s} \frac{x}{s^\frac{1}{\alpha}} ds \leq C x.
\end{align*}
and thus a bound on the behaviour of $k$ at $0$:
\begin{align*}
    \text{if } \alpha \leq 1, \quad  \vert k(x) \vert \leq C x^{\alpha-1} \vert \ln(x) \vert; \quad \text{if } \alpha \geq 1, \quad  \vert k(x) \vert \leq C. 
\end{align*}
This last inequality justifies $k\in L^1(\R)$, and justifies that $k$ is bounded when $x$ goes to $0$ for $\alpha>1$.

In addition for $\alpha>1$, we obtain more results on $k'$. We prove first the asymptotic development at $+\infty$, and its proof is along the same lines as the arguments for $k$. We have for any $x\neq 0$, by the development \eqref{estimation:h_prime} for the integrability at $0$:
\begin{align*}
    k'(x)=\frac{1}{\pi}\int_{0}^{+\infty}\frac{e^{-s}}{s^{\frac{2}{\alpha}}}h'\left( \frac{x}{s^{\frac{1}{\alpha}}} \right) ds .
\end{align*}

Performing as for $k$, we get the asymptotic development of $k'$ for $\alpha>1$. We investigate on the behaviour of $k'$ at $0$. As before, we split the integral into two parts: on $(0,x^\alpha)$, by \eqref{estimation:h_prime}, the first part of the integral is bounded by $x^{\alpha-2}$; on $(x^\alpha,+\infty)$, with \eqref{eq:bound_h}, this second part of the integral is bounded by $x^{\alpha-2}$. With the two estimates, $k'$ is in $L^1(\R)$ for $\alpha>1$.

This concludes the proof of Lemma \ref{lemma:est:invlapl}.

\end{proof}

\begin{coro}\label{lemma:est:kernel}
Let $g\in L^1(\R)$ satisfying $\vert g(x) \vert \leq C | x |^{-\alpha-1}$. There exists $C=C(g)$ such that: 
\begin{align}\label{est:kernel}
    |k\ast g|(x)\leq \frac{C}{\langle x \rangle^{\alpha+1}}.
\end{align}

Furthermore, if $g\in \mathcal{C}^1(\R)$ and $\vert g'(x) \vert \leq C | x |^{-2-\alpha}$, then there exists $C=C(g,g')$ such that: 
\begin{align}
    \vert \partial_x \left( k \star g\right)\vert (x) \leq \frac{C}{ \langle x\rangle ^{\alpha+2}}.
\end{align}
\end{coro}

\begin{proof}
For sake of simplicity, we focus on the case $x\geq 0$, the negative terms are obtained by parity of $k$. First, by the decomposition of \eqref{defi:omega_x}, we have that:
\begin{align*}
    k\ast g(x)= \int_{\Omega_x} k(x-y)g(y)dy+\int_{\Omega_x^c} k(x-y)g(y)dy=I_1+I_2. 
\end{align*}
By the decay assumption on $g$, we get that:
\begin{align*}
    |I_2|\leq \frac{C}{\langle x\rangle^{\alpha+1}}.
\end{align*}
Moreover, from Lemma \ref{lemma:est:invlapl}, we deduce that:
\begin{align*}
    |k(x)|\leq \frac{C}{\langle x\rangle^{\alpha+1}}.
\end{align*}
Furthermore, the inequality $|y|\leq \displaystyle\frac{x}{2}$ implies $|x-y|\geq \displaystyle\frac{x}{2}$. Thus, we deduce that: 
\begin{align*}
    |I_1|\leq \frac{C}{\langle x\rangle^{\alpha+1}}.
\end{align*}
Gathering these estimates, we conclude \eqref{est:kernel} for $x>0$. By arguing similarly as the case $x>0$, we get \eqref{est:kernel} for $x\in \R$.

The second estimate is based on one step further on the asymptotic development of $k$:
\begin{align}\label{est:kgrogro}
    \left\vert k(x) - \frac{k_1}{x^{\alpha+1}} \right\vert \leq \frac{C}{\langle x\rangle^{2\alpha+1}}.
\end{align}
By the same decomposition as previously, we have:
\begin{align}
    \partial_x (k \star g)(x)= I_1 +I_2,
\end{align}
where: 
\begin{align*}
    I_1 & = \int_{\Omega_x} \left(  k(x-y) - \frac{k_1}{(x-y)^{\alpha+1}} \right) g'(y) dy + \int_{\Omega_x} \frac{k_1}{(x-y)^{\alpha+1}} g'(y) dy. \\
    I_2 & = \int_{\Omega_x^c} k(x-y) g'(y) dy.
\end{align*}
It suffices to prove the bound with $I_1$, since $I_2$ is dealt with like the previous step. For the first term, due to \eqref{est:kgrogro}, we have:
\begin{align*}
    \left\vert \int_{\Omega_x} \left(  k(x-y) - \frac{k_1}{\vert x- y \vert^{\alpha+1}} \right) g'(y) dy \right\vert \leq \frac{C}{ \langle x\rangle^{2\alpha+1}}.
\end{align*}
For the second term of $I_1$, by integration by part:
\begin{align}
    \left\vert \int_{\Omega_x} \frac{k_1}{(x-y)^{\alpha+1}} g'(y) dy \right\vert \leq \left\vert \frac{k_1}{(\frac{x}{2})^{\alpha+1}} g\left(-\frac{x}{2}\right)\right\vert +\left\vert \frac{k_1}{(\frac{3x}{2})^{\alpha+1}} g\left(\frac{x}{2}\right)\right\vert + \left\vert \int_{\Omega_x} \frac{k_1(\alpha+1)}{(x-y)^{\alpha+2}} g(y) dy \right\vert \leq \frac{C}{\langle x\rangle^{\alpha+2}}.
\end{align}

This achieves the proof of one derivative applied to $k\star g$ in Lemma \ref{lemma:est:kernel}.
\end{proof}

\section{Proof of the asymptotic expansion of order 1}

In this section we prove Theorem \ref{propo:notre_method_1}, the asymptotic development of order 1 of $Q$ for $\alpha\in(0,2)$ and $p\in (1,p^*(\alpha))$. \cite{FLS16} proved the decay of a positive $Q$ at infinity:
\begin{align}\label{estimation:Q}
    |Q(x)|\leq \frac{C}{\langle x\rangle ^{\alpha+1}}.
\end{align}
We propose here to establish a more precise asymptotic development in the general case that $Q$ is not necessarily positive, but satisfies \eqref{hypo:Q}. The proof is separated in two steps: the first consists in establishing the first order asymptotic, and to prove in particular that \eqref{estimation:Q} is satisfied. Second, we get successively the asymptotic development of order $1$ of the derivatives of $Q^{(j)}$, by using an induction process on $j$.

\subsection{First order expansion of $Q$}\label{subsec:first_order_Q}
Let us prove the asymptotic expansion of order 1 for $ \vert x \vert >1$. To begin with, if $Q$ is a weak solution of \eqref{defi:eq_elliptic}, then $Q$ is a weak solution of the following equation:
\begin{align}\label{definition:Q_p}
    Q= \mathcal{F}^{-1}\left( \left( \vert \xi \vert^\alpha + 1 \right)^{-1} \right) \star \left(|Q|^{p-1}Q\right) =k \star \left(|Q|^{p-1}Q\right).
\end{align}
Before giving the asymptotic expansion, we improve the bound on $Q$:
$$
|Q(x)|\leq\frac{C}{\langle x \rangle^{\alpha+1}}.
$$
 We split the integral \eqref{definition:Q_p} in high and low values, with $\Omega_x$ defined in \eqref{defi:omega_x}:
\begin{align}\label{defi:K1_K2}
    K_1:= \int_{ \Omega_x} k(x-y) |Q|^{p-1}Q(y)dy, \quad K_2 := \int_{ \Omega_x^c} k(x-y)|Q|^{p-1}Q(y) dy.
\end{align}
Using \eqref{est:invlapl} on $K_1$, and the hypothesis $Q\in L^p(\R)$, we deduce that:
\begin{align*}
    |K_1|\leq \frac{C}{\langle x\rangle^{\alpha+1}}.
\end{align*}
Now we continue with $K_2$. We use the decay assumption on $Q$: $|Q(x)|\leq C \langle x \rangle^{-l}$, and $k\in L^1(\R)$ to deduce:
\begin{align*}
    |K_2|\leq \frac{C}{\langle x \rangle ^{pl}}.
\end{align*}
Then, we obtain that $|Q(x)|\leq C \langle x \rangle ^{-pl}$. Since $p>1$, by iterating the previous steps with the improved bound on the decay of $Q$, we conclude that:
\begin{align*}
    |Q(x)|\leq \frac{C}{\langle x \rangle^{\alpha+1}}
\end{align*}

From this inequality, we deduce that $Q\in L^q(\R)$ for any $q\in [1,+\infty]$. By the equation \eqref{defi:eq_elliptic}, $\vert D\vert^\alpha Q$ is also in $L^q(\R)$ and $Q\in H^\alpha(\R)$. Thus $Q$ is a solution of \eqref{defi:eq_elliptic} and solves \eqref{definition:Q_p} in the usual sense.

We continue with the asymptotic expansion of $Q$. We find the equivalent at the first order, which comes from $K_1$. We have:
\begin{align}
    \left\vert K_1 - \int_{\Omega_x} \frac{k_1}{(x-y)^{1+\alpha}} \vert Q\vert^{p-1}Q(y) dy \right\vert & \leq \int_{\Omega_x} \left\vert k(x-y) - \frac{k_1}{(x-y)^{\alpha+1}} \right\vert \vert Q \vert ^{p}(y)dy \nonumber\\
        & \leq C \int_{\Omega_x} \frac{\vert Q\vert ^{p}(y)}{\vert x-y \vert^{2\alpha+1}}dy \leq \frac{C}{ \langle x \rangle ^{2\alpha+1}}. \label{eq:K1_1}
\end{align}
We then get a development in $x$ of the term in the left hand-side, where $a_1$ is defined in \eqref{defi:a1_a2_a3}:

\begin{align*}
\left\vert \int_{ \Omega_x} \frac{k_1}{(x-y)^{\alpha+1}} \vert Q \vert^{p-1} Q (y) dy - \frac{a_1}{ x^{\alpha+1}} \right\vert \leq \frac{C}{\langle x\rangle^{\alpha+1}} \int_{\Omega_x^c} \vert Q \vert^p(y) dy + \frac{C}{\langle x\rangle ^{\alpha+2}} \int_{\Omega_x} \vert y \vert \vert Q \vert^{p}(y) dy.
\end{align*}
The first term is dealt with the asymptotic bound $Q(y) \leq C \langle y \rangle^{-(1+\alpha)}$:
\begin{align*}
    \frac{1}{\langle x \rangle ^{\alpha+1}} \int_{\Omega_x^c} \vert Q\vert ^{p}(y) dy \leq \frac{C}{\langle x \rangle ^{\alpha+p(1+\alpha)}}.
\end{align*}
For the second term, we need to split the different cases between $\alpha$ and $p$:
\begin{align*}
    \frac{C}{\langle x \rangle ^{\alpha+2}} \int_{\Omega_x} \vert y \vert^{1-\min(\frac{p\alpha}{2},1)} \vert y \vert ^{\min(\frac{p\alpha}{2},1)} \vert Q \vert^{p}(y) dy \leq \frac{C}{\langle x \rangle ^{\alpha+1 +\min(\frac{p\alpha}{2},1 )}}.
\end{align*}

We continue with $K_2$:
\begin{align}\label{eq:K2}
    \left\vert K_2 \right\vert \leq C \int_{ \Omega_x^c} \frac{\vert k(x-y) \vert}{\langle x \rangle ^{p(\alpha+1)}} dy \leq \frac{C}{\langle x \rangle ^{p(\alpha+1)}}.
\end{align}
Gathering the estimates on $K_1$ and on $K_2$ we obtain the asymptotic development of $Q$ at the main order:
\begin{align}
    \forall x >1, \quad \left\vert Q(x)- \frac{a_1}{  x  ^{\alpha+1}} \right\vert \leq \frac{C}{\langle x\rangle ^{2\alpha+1}} + \frac{C}{\langle x \rangle ^{p(\alpha+1)}}+\frac{C}{\langle x\rangle^{\alpha+1 +\min(\frac{p\alpha}{2},1 )}}.
\end{align} 

\subsection{First order expansion of the derivatives of $Q$}

This section is dedicated to the asymptotic expansion of any $j$-derivative of $Q$, for $\alpha>1$ and $j \leq \lfloor p \rfloor$. We use the notation $f(u):=\vert u \vert^{p-1}u$. We prove the following statement by induction on $j\in [1, \lfloor p \rfloor]$:
\begin{stat}\label{stat:induction_Q}
    The function $Q$ is of class $\mathcal{C}^j(\R)$. If $j\geq 2$, there exists a sequence of polynomials $R_{j,l}$ such that:
    \begin{equation}\label{eq:R_jk}
        \frac{d^j}{dx^j}(f(Q)) = p Q^{(j)} f'(Q)+ \sum_{l=2}^j R_{j,l}(Q',\cdots, Q^{(j-1)}) f^{(l)} (Q),
    \end{equation}
    where $R_{j,l}$ satisfies:
    \begin{align*}
        \forall j,l, \quad \vert R_{j,l}(Q'(x), \cdots, Q^{(j-1)}(x)) f^{(l)}(Q(x)) \vert \leq \frac{C}{\langle x \rangle^{p(\alpha+1)+j}}.
    \end{align*}
    If $j\geq 1$, the following asymptotic expansion holds:
    \begin{align}\label{eq:decay_Q_j}
    \exists \beta>\alpha+1+j, \forall x \geq 1, \quad  Q^{(j)}(x) -(-1)^j\frac{(\alpha+j)!}{\alpha!} \frac{a_1}{x^{\alpha+1+j}}  =o_{+\infty}\left(x^{-\beta}\right).
    \end{align}
\end{stat}

Notice that in this induction process, we need $j$ such that $p-j>0$ to ensure that $f^{(p-j)}(Q)$ or equivalently $\vert Q\vert ^{p-j}$ is well-defined; this constraint does not occur if $Q$ is positive.

We prove the statement for $j=1$. The first derivative of $Q$ is continuous since $k'$ is in $L^1(\mathbb{R})$ (see Claim \ref{claim:k_function_of_h}), $Q\in L^q(\R)$, for any $q\in [1,+\infty]$ (see Theorem \ref{propo:notre_method_1}), and given by $ Q'= k' \star (|Q|^{p-1}Q) $. Moreover, we have $Q'\in L^{\infty}(\R)$, since $k'\in L^1(\R)$ and $Q\in L^{\infty}(\R)$. Now, we prove that:
\begin{align*}
    |Q'(x)|\leq \frac{C}{\langle x \rangle^{\alpha+2}}.
\end{align*}
We write, for $x>1$:
\begin{align}\label{defi:J1_J2}
    J_1:=\int_{\Omega_x} k'(x-y)|Q(y)|^{p-1}Q(y)dy, \quad    J_2:=\int_{\Omega_x^c} k'(x-y)|Q(y)|^{p-1}Q(y)dy. 
\end{align}
On the one hand, by using $Q\in L^{p}(\R)$ and the decay property of $k'$ \eqref{est:invlapl_deriv}, we deduce that: 
\begin{align*}
    |J_1|\leq \frac{C}{\langle x\rangle^{\alpha+2}}.
\end{align*}
On the other hand, from the bound $|Q(x)|\leq C \langle x\rangle ^{-(1+\alpha)}$ and the fact $k'\in L^1(\R)$ for $\alpha>1$, we have that:
\begin{align*}
    |J_2|\leq \frac{C}{\langle x \rangle^{p(1+\alpha)}}. 
\end{align*}
Then, we obtain that $|Q'(x)|\leq C\langle x\rangle^{-\min(p(1+\alpha), 2+\alpha)}$. By integrating by parts $J_2$, and using the new bound on $Q'$, we get that:
\begin{align*}
    |Q'(x)|\leq \frac{C}{\langle x \rangle^{\min((2p-1)(1+\alpha),2+\alpha)}}.
\end{align*}
Since $p>1$, by iterating the previous steps, we conclude that: 
\begin{align*}
    |Q'(x)|\leq \frac{C}{\langle x \rangle^{2+\alpha}}.
\end{align*}
Now, we derive the asymptotic development of order 1 of $Q'$ in \eqref{eq:decay_Q_j}, with:
\begin{align*}
    Q'(x) = k\star ((|Q|^{p-1}Q)')(x).
\end{align*}
First, we estimate the convolution on $\Omega_x$. From the asymptotic development of $k$ \eqref{est:invlapl}, and the fact $Q'|Q|^{p-1}\in L^1(\R)$ we obtain that:
\begin{align*}
    \left|\int_{\Omega_x}\left(k(x-y) - \frac{ k_1}{(x-y)^{\alpha+1}} - \frac{ k_2}{(x-y)^{2\alpha+1}}\right)(|Q|^{p-1}Q)'(y) dy \right|\leq \frac{C}{\langle x\rangle ^{3\alpha+1}}.
\end{align*}
By integrating by parts, we deduce that:
\begin{align*}
    \left|\int_{\Omega_x} \frac{ k_2}{(x-y)^{2\alpha+1}}(|Q|^{p-1} Q)'(y) dy \right|\leq \frac{C}{\langle x\rangle^{2\alpha+2}}. 
\end{align*}
After an integration by parts on the term with $k_1$, using the fact $\vert y\vert ^{\min\left(1,\frac{p\alpha}{2}\right)}Q\in L^{1}(\R)$, we get that:
\begin{align*}
    \left|\int_{\Omega_x}\frac{k_1}{(x-y)^{\alpha+1}}(|Q|^{p-1}Q)'(y)dy +\frac{(\alpha+1)k_1}{x^{\alpha+2}}\int_{\Omega_x}|Q|^{p-1}(y) Q(y)dy \right|\leq \frac{C}{\langle x \rangle^{\alpha+2+\min\left(1,\frac{p\alpha}{2}\right)}} + \frac{C}{\langle x \rangle^{2\alpha+2}}.
\end{align*}
Furthermore, using the decay assumption on $Q$, we obtain that:
\begin{align*}
    \left|\frac{(\alpha+1)k_1}{x^{\alpha+2}}\int_{\Omega_x^c}|Q|^{p-1}(y) Q(y)dy \right|\leq \frac{C}{\langle x\rangle^{(p+1)(\alpha+1)+1}}.
\end{align*}
To summarize, we have proved that:
\begin{align*}
    \exists \beta>\alpha+2, \quad \int_{\Omega_x} k(x-y)(|Q|^{p-1}(y)Q(y))' dy + \frac{(\alpha+1)k_1}{x^{\alpha+2}}\int_{\R}|Q|^{p-1}(y)Q(y)dy = o_{+\infty}\left(x^{-\beta} \right).
\end{align*}
To finish the proof of the asymptotic development of $Q'$, we have to estimate $Q'(x) = k\star ((|Q|^{p-1}Q)')(x)$ on $\Omega_x^c$. By using the decay assumption on $Q$ and $Q'$, we obtain that:
\begin{align*}
    \left|\int_{\Omega_x^c} k(x-y)|Q|^{p-1}(y)Q'(y)dy \right|\leq \frac{C}{\langle x\rangle^{p(1+\alpha)+1}}.
\end{align*}
Therefore, gathering these estimates, we conclude that:
\begin{align*}
    \exists \beta>\alpha+2, \quad Q'(x) + \frac{(\alpha+1)a_1}{x^{\alpha+2}}=o_{+\infty}\left(x^{-\beta} \right).
\end{align*}

To prove the statement for the other $j$, we proceed by induction. The regularity of $Q^{(j+1)}$ is obtained by $k'\in L^1(\mathbb{R})$, the formula \eqref{eq:R_jk} at $j$, and the asymptotic expansion of $Q^{(i)}$ for $i \leq j$. Thus, we have:
\begin{align*}
    Q^{(j+1)}= k \star (\vert Q\vert^{p-1}Q)^{(j+1)}.
\end{align*}

The proof of the existence of the polynomials $R_{j+1,l}$ in \eqref{eq:R_jk}, for $l\in [2,j+1]$, is given by direct computations. Before giving the bound on the polynomials, we remark that for all $l\in\{0,\cdots, \lfloor p \rfloor\}$ there exists $C_l>0$ such that:
\begin{align*}
    f^{(l)}(x)=\begin{cases}C_lx|x|^{p-l-1}, \quad  &\text{ if $l$ is even},\\
    C_l|x|^{p-l}, \quad  &\text{ if $l$ is odd }
    \end{cases}.
\end{align*}

We detail the bound satisfied by the polynomial $R_{j+1,l}$. If $l=2$, we have that:
\begin{align*}
    |Q^{(j)}f^{(2)}(Q)| \leq \frac{C}{\langle x \rangle^{p(1+\alpha)+j+1}} \quad \text{and} \quad \left\vert \frac{d}{dx}(R_{j,2})f^{(2)}(Q)\right\vert \leq \frac{C}{\langle x \rangle^{p(1+\alpha)+j+1}}.
\end{align*}
Moreover, for a fixed number $l\in [3,j+1]$, the polynomial $R_{j+1,l}$ is at most composed of two terms:
\begin{align*}
    \left\vert R_{j,l-1} f^{(l)}(Q) \right\vert \leq \frac{C}{\langle x \rangle^{p(1+\alpha)+j+1}} \quad \text{and} \quad  \left\vert \frac{d}{dx}(R_{j,l})f^{(l)}(Q)\right\vert \leq \frac{C}{\langle x \rangle^{p(1+\alpha)+j+1}},
\end{align*}

which concludes the decay of the term with $f^{(l)}$ for $j+1$.

Now, we prove \eqref{eq:decay_Q_j} for $Q^{(j+1)}$. Let us sketch the proof of the rough bound of $Q^{(j+1)}$:
\begin{align}\label{eq:bound_Q_j+1}
    \vert Q^{(j+1)}(x) \vert \leq \frac{C}{\langle x \rangle^{\alpha+2+j}}.
\end{align}
As for the first derivative, we define $J_{1,j+1}$ and $J_{2,j+1}$ as in \eqref{defi:J1_J2}:
\begin{align}
    J_{1,j+1}:=\int_{\Omega_x} k'(x-y)(|Q|^{p-1}Q)^{(j)}(y) dy, \quad    J_{2,j+1}:=\int_{\Omega_x^c} k'(x-y)(|Q|^{p-1}Q)^{(j)}(y)dy. 
\end{align}

For the first term with the asymptotic development \eqref{est:invlapl_deriv} of $k'$, for any $N$ and by successive integration by parts:
\begin{align*}
    \left\vert J_{1,j+1} \right\vert & \leq \left\vert \int_{\Omega_x}\left( k'(x-y) -\sum_{n=1}^N \frac{(n\alpha+1)k_n}{(x-y)^{n\alpha+2}} \right) (|Q|^{p-1}Q)^{(j)}(x) dx \right\vert \\
        & \quad + \left\vert \int_{\Omega_x} \sum_{n=1}^N \frac{(n\alpha+1)k_n}{(x-y)^{n\alpha+2}} (|Q|^{p-1}Q)^{(j)}(x) dx \right\vert \\
    & \leq \frac{C}{\langle x\rangle ^{(N+1)\alpha+2}} + \frac{C}{\langle x\rangle ^{\alpha+2+p(\alpha+1)+j}} + \frac{C}{\langle x \rangle^{\alpha+2+j}} \leq \frac{C}{\langle x \rangle^{\alpha+2+j}},
\end{align*}
where the last inequality holds for $N$ large enough, for instance $N>\frac{j}{\alpha}$. 

For the second term, by Lemma \ref{lemma:est:invlapl} and the Statement, we have:
\begin{align*}
    \left\vert J_{2,j+1} \right\vert \leq \frac{C}{\langle x \rangle ^{p(\alpha+1)+j}}.
\end{align*}

Once we have the bound $\vert Q^{(j+1)}(x)\vert \leq C\langle x \rangle ^{-(\alpha+2+j)} +C\langle x\rangle ^{-(p(\alpha+1)+j)}$, we improve this bound by injecting it successively in $J_{2,j+1}$ after an integration by part and with \eqref{eq:R_jk}:
\begin{align*}
    \left\vert J_{2,j+2} \right\vert \leq C \| ( \vert Q \vert ^{p-1} Q)^{(j+1)} \|_{L^\infty(\Omega_x^C)} \leq \frac{C}{\langle x \rangle ^{(p-1)(\alpha+1)}} \| Q^{(j+1)} \|_{L^\infty(\Omega_x^C)} + \frac{C}{\langle x \rangle ^{p(\alpha+1)+j+1}}.
\end{align*}

By doing the injection enough times, we get the bound \eqref{eq:bound_Q_j+1}.

To obtain the exact asymptotic development, it suffices to use the same proof as the one of $Q'$ and the decomposition used in $J_{1,j+1}$ to extract the term of main order. We get:
\begin{align*}
    \exists \beta=\beta(j)>\alpha+j+2, \quad \left\vert J_{1, j+1} - \frac{k_1}{x^{\alpha+j+2}}\frac{(\alpha+j+1)!}{\alpha!} \int \vert Q \vert^p Q \right\vert =o_{+\infty} \left(x^{-\beta}\right)
\end{align*}
By injecting \eqref{eq:bound_Q_j+1} in $J_{2,j+1}$, we conclude the asymptotic expansion \eqref{eq:decay_Q_j} for $j+1$.

By finishing to prove the Statement for any $j\leq \lfloor p\rfloor$, we conclude the proof of Theorem \ref{propo:notre_method_1}.

Notice that by the same method, we obtain the existence of the next derivative:
\begin{align}
    Q^{(\lfloor p\rfloor+1)} = \int k'(x-y) \left( \vert Q \vert^{p-1} Q\right)^{(\lfloor p \rfloor)} (y) dy.
\end{align}
However, an integration by parts is not justified since $Q$ can be equal to $0$ at some points, therefore the function $\left( \vert Q \vert^{p-1} Q\right)^{(\lfloor p \rfloor+1)}$ may not be defined. Our method does not provide an asymptotic expansion of $Q^{(\lfloor p\rfloor +1)}$.

\section{Asymptotic expansion for positive solutions} 

To start this section dedicated to the proof of Theorem \ref{propo:notre_method_2}, we recall that, on the one hand, if $Q>0$ then $Q\in H^{\infty}(\R)$, see \cite{FLS16}, the third point in the Remarks after Proposition 1.1, and $Q$ is even (up to translation) see Proposition 1.1 of \cite{FLS16} and \cite{MZ10}. 

Since $Q\in H^{\infty}(\R)$ and even, we give the asymptotic expansion of $Q^{(j)}$ for any $j$. The proof of this result follows the line of the one of Theorem \ref{propo:notre_method_1}. We use the induction on the statements \eqref{eq:R_jk} and \eqref{eq:decay_Q_j} as for Theorem \ref{propo:notre_method_1}. For $Q>0$ the induction process given in Theorem \ref{propo:notre_method_1} does not stop anymore. Therefore, in the case of $Q>0$ we get the asymptotic development of order 1 for all the derivatives.

We continue by giving the second order expansion of $Q$.

The proof is based on the same arguments as the ones of subsection \ref{subsec:first_order_Q}, where we have obtained:
\begin{align*}
    \exists \beta> \alpha+1, \quad Q(x)- \frac{a_1}{x^{\alpha+1}} = o_{+\infty}\left(x^{-\beta} \right).
\end{align*}
We prove in this section the expansion at the next order. We study separately the different cases.

Let $K_1$ and $K_2$ be defined as in \eqref{defi:K1_K2}, we have the following decomposition. First, with $K_1$:
\begin{align*}
    K_1 =K_{11} + K_{12} + K_{13} + K_{14}
\end{align*}
with
\begin{align*}
      K_{11} & = \int_{\Omega_x} \left( k(x-y) - \frac{k_1}{(x-y)^{\alpha+1}} - \frac{k_2}{(x-y)^{2\alpha+1}} \right) Q^p(y)dy, \\
      K_{12} & = \int_{\Omega_x} \left( \frac{k_1}{(x-y)^{\alpha+1}} + \frac{k_2}{(x-y)^{2\alpha+1}} - \left( \frac{k_1}{x^{\alpha+1}} + \frac{k_2}{x^{2\alpha+1}} \right) \right) Q^p(y)dy, \\
      K_{13} & =\left( \frac{k_1}{x^{\alpha+1}} +\frac{k_2}{x^{2\alpha+1}} \right) \int_\mathbb{R} Q^p(y) dy \quad \text{and} \quad K_{14} = - \left( \frac{k_1}{x^{\alpha+1}} +\frac{k_2}{x^{2\alpha+1}} \right) \int_{\Omega_x^c} Q^p(y) dy.
\end{align*}

By the asymptotic expansion \eqref{est:invlapl} of $k$, we have $\vert K_{11} \vert \leq C\langle x\rangle^{-(3\alpha+1)}$. For $K_{12}$, we use the parity of $Q$ to get $\int_{\Omega_x} y Q^p(y)dy=0$, therefore we obtain:
\begin{align*}
    \left\vert K_{12} \right\vert \leq \frac{C}{\langle x \rangle^{3+\alpha}} \int y^2 Q^p(y)dy + \frac{C}{\langle x \rangle^{3+\alpha}} \int \vert y \vert Q^p(y) dy \leq \Theta(x)
\end{align*}
where
\begin{equation*}
    \Theta(x) =
    \left\{ \begin{aligned}
        & \frac{1}{\langle x\rangle ^{\alpha+3}} & \quad \text{if} \quad p(1+\alpha)>3,\\
        & \frac{\ln(\vert x\vert )}{\langle x\rangle^{\alpha+3}} & \quad \text{if} \quad p(1+\alpha)=3,\\
        & \frac{1}{\langle x \rangle^{\alpha+p(1+\alpha)}} & \quad \text{if} \quad p(1+\alpha)<3. \\
    \end{aligned} \right.
\end{equation*}
From the asymptotic of $Q$ on $\Omega_x^c$, we have:
\begin{align*}
    \vert K_{14} \vert \leq \frac{C}{x^{p(1+\alpha)+\alpha}}.
\end{align*}

$K_{13}$ is the only remaining term that could potentially give the next order term in the asymptotic expansion of $Q$.

We decompose $K_2$ as:
\begin{align*}
    K_2=\int_{\Omega_x^c\cap (\frac{x}{2},\frac{3x}{2})} k(x-y)Q^{p}(y)dy + \int_{\Omega_x^c\cap (\frac{x}{2},\frac{3x}{2})^{c}} k(x-y)Q^p(y)dy=K_{2,1}+K_{2,2}.
\end{align*}
Using the decay assumptions on $k$ \eqref{est:invlapl} and $Q$ given in Theorem \ref{propo:notre_method_1} we obtain that:
\begin{align*}
    |K_{2,2}|=\left|\int_{\Omega_x^c\cap (\frac{x}{2},\frac{3x}{2})} k(x-y)Q^{p}(y)dy \right|\leq \frac{C}{\langle x \rangle^{p(1+\alpha)+\alpha}}.
\end{align*}
Furthermore, by using again the decay estimate on $k$ \eqref{est:invlapl} and the asymptotic expansion of $Q$ given in Theorem \ref{propo:notre_method_1}, we deduce that for some $\beta>p(1+\alpha)$:

\begin{align*}
    K_{2,1}= \int_{\Omega_x^c\cap (\frac{x}{2},\frac{3x}{2})} k(x-y)\frac{a_1^p}{y^{p(\alpha+1)}}\left(1+\frac{Q(y)- \frac{a_1}{y^{\alpha+1}}}{\frac{a_1}{y^{\alpha+1}}}\right)^{p}dy = \frac{\tilde{a}_1}{x^{p(1+\alpha)}}+ o_{+\infty}\left(x^{-\beta}\right),
\end{align*}
with $\tilde{a}_1:= a_1^p\displaystyle\int_\R k(x)dx$, and let us decompose the proof into three cases, depending on $p$ with respect to $\frac{2\alpha+1}{\alpha+1}$. 

\textbf{First case : $p<\frac{2\alpha+1}{\alpha+1}$}

This case corresponds to a low non-linearity compared to the influence of the dispersion, and the biggest error term comes from $K_{21}$. Gathering the estimates on $K_1$ and on $K_2$, we obtain:
\begin{align*}
    \exists \beta > p(\alpha+1), \quad Q(x) -\frac{a_1}{x^{\alpha+1}} -\frac{\tilde{a}_1}{x^{p(\alpha+1)}}  = o_{+\infty} \left( x^{-\beta} \right)
\end{align*}

\textbf{Second case: $p=\frac{2\alpha+1}{\alpha+1}$}

In this particular case of balance between dispersion and non-linearity, the next order is given by two different terms, from $K_{13}$ and $K_{21}$:
\begin{align*}
    \exists \beta > 2\alpha+1 , \quad Q(x)- \frac{a_1}{x^{\alpha+1}} - \frac{\tilde{a}_1}{x^{2\alpha+1}} - \frac{a_2}{x^{2\alpha+1}} =o_{+\infty}\left( x^{-\beta}\right).
\end{align*}
where $a_2$ is given in \eqref{defi:a1_a2_a3}.

\textbf{Third case: $p>\frac{2\alpha+1}{\alpha+1}$} 

When the non-linearity is above $\frac{2\alpha+1}{\alpha+1}$, the tail of $Q$ is negligible compared to the next order term given by the dispersion:
\begin{align*}
     \exists \beta > 2\alpha+1, \quad Q(x)- \frac{a_1}{x^{\alpha+1}} - \frac{a_2}{x^{2\alpha+1}} =  o_{+\infty}\left( x^{-\beta} \right).
\end{align*}

 This concludes the proof of Theorem \ref{propo:notre_method_2}.

\section{Asymptotic expansion for polynomial non-linearities}

The proof of Theorem \ref{propo:poly_derivatives} follows the arguments given in the proof of the Theorem \ref{propo:notre_method_1} and of Theorem \ref{propo:notre_method_2}.

We now continue with the proof of Theorem \ref{thm:asympQ}. In this section the non-linearity is fixed at $p=3$ and the dispersion $\alpha\in (1,2)$.

 To begin with, using \eqref{definition:Q_p}, we obtain the decomposition in high and low values:  
\begin{align*}
    Q(x)=\int_{\Omega_x} k\left(x-y\right) Q^{3}\left(y\right) dy + \int_{\Omega_x^c} k\left(x-y\right)Q^{3}\left(y\right) dy=: K_1 + K_2 .
\end{align*}
Using the asymptotic development of $k$ \eqref{est:invlapl}, we get that:
\begin{align}
\MoveEqLeft
    \left\vert K_1 - \int_{\Omega_x}\left( \frac{k_1}{ (x-y)^{\alpha+1}}  +\frac{k_2}{ (x-y)^{2\alpha+1}} \right) Q^3(y)dy\right\vert \notag \\
        & = \left|\int_{\Omega_x} \left( k(x-y)- \frac{k_1}{ (x-y)^{\alpha+1}} - \frac{k_2}{ (x-y)^{2\alpha+1}} \right) Q^3(y)dy \right| \notag \\
        & \leq C \int_{\Omega_x} \frac{1}{|x-y|^{3\alpha+1}} Q^3(y) dy \leq \frac{C}{\langle x \rangle ^{3\alpha+1}}. \label{eq:zfraction_DL_x}
\end{align}
From the asymptotic expansion of $\frac{1}{(1+x)^{\beta}}$ for $\beta\in\mathbb{R}$, with $\int \vert y \vert ^3 Q^3(y) dy<\infty$ for $\alpha>1$, we deduce that:
\begin{align}
\MoveEqLeft
    \left| \int_{\Omega_x}\left( \frac{k_1}{ (x-y)^{\alpha+1}} -k_1 \left( \frac{1}{x^{\alpha+1}} + (\alpha+1)\frac{y}{x^{\alpha+2}}+\frac{(\alpha+1)(\alpha+2)}{2} \frac{y^2}{x^{\alpha+3}} + \frac{(\alpha+1)(\alpha+2)(\alpha+3)}{3!} \frac{y^3}{x^{\alpha+4}} \right) \right. \right. \notag \\
        & \quad \quad \left. \left.+ \frac{k_2}{ (x-y)^{2\alpha+1}} - k_2 \left( \frac{1}{x^{2\alpha+1}}  +(2\alpha+1)\frac{y}{x^{2\alpha+2}}\right) \right)Q^3(y)dy \right|\notag \\
        & \quad \quad  \leq \frac{C}{\langle x\rangle^{\alpha+5}} + \frac{C}{\langle x\rangle^{2\alpha+3}} \leq \frac{C}{\langle x\rangle^{2\alpha+3}}  .\label{asympt:2}
\end{align}

Since $\displaystyle\int_{\Omega_x}yQ^{3}(y)dy=\int_{\Omega_x}y^3 Q^3(y)dy=0$, we deduce that the terms for the asymptotic development of $Q$ are  $\displaystyle \frac{a_1}{x^{\alpha+1}} + \frac{a_2}{x^{2\alpha+1}} + \frac{a_3}{x^{\alpha+3}}$, with $a_1,a_2,a_3$ defined in \eqref{defi:a1_a2_a3}. Now, we have to verify that the estimates on the asymptotic development hold on $\Omega_x^c$.

From the decay of $Q$ \eqref{estimation:Q}, we obtain that: 
\begin{align}
    |K_2|\leq \int_{\Omega_x^c} \vert k(x-y)\vert dy \| Q \|_{L^{\infty}(\vert y\vert \geq \frac{x}{2})}^3 \leq \frac{C}{ \langle x \rangle  ^{3(1+\alpha)}} \quad 
\text{ and }\quad 
    \int_{\Omega_x^c} \frac{ Q^3(y)}{x^{3(\alpha+1)}} + \frac{y^2Q^3(y)}{x ^{\alpha+3}}  dy \leq \frac{C}{\langle x \rangle ^{4\alpha+3}}.
\end{align}
This concludes the proof of the estimate \eqref{asympt:Q}.

The proof of the estimate \eqref{asympt:Q'} is similar as the proof of the estimate \eqref{asympt:Q}.

\appendix 

\section{Regularity result for polynomial non-linearities}\label{appendix_A}

Let $p\in\N$. We prove in this appendix that if $Q$ verifies \eqref{hypo:Q}, then:
\begin{align*}
    \forall \beta \in \mathbb{R}^+, \quad \| Q\|_{H^{\beta}} <\infty.
\end{align*}
We prove this statement by induction. Since $Q$ is in $L^q(\R)$ for any $q\in [1,+\infty]$ by Theorem \ref{propo:notre_method_1}, we obtain the result for $\beta=0$. To prove that $Q \in H^{(n+1)\alpha}(\R)$ with the assumption $Q \in H^{n\alpha}(\R)$, it suffices to study $\vert D \vert^{n\alpha} (Q^p)$. By the fractional Leibniz rule (also called Kato-Ponce commutator estimate, see \cite{GO14} for the endpoint), and $\gamma>0$, we have:
\begin{align*}
    \| \vert D \vert ^\gamma (Q^p) \|_{L^2} 
        & \leq \| \langle \vert D \vert \rangle^\gamma (Q^p) \|_{L^2} \leq C \left( \| Q \|_{L^\infty} \| \langle \vert D \vert \rangle ^{\gamma}(Q^{p-1}) \|_{L^2} + \| \langle \vert D \vert \rangle^{\gamma} Q \|_{L^2} \| Q^{p-1} \|_{L^\infty}  \right) \\
        & \leq C \| Q \|_{L^\infty}^{p-1} \| \langle \vert D \vert \rangle^{\gamma}  Q \|_{L^2},
\end{align*}
where the last step is obtained by induction. Thus we obtain:
\begin{align*}
    \| \vert D \vert^{(n+1)\alpha} Q \|_{L^2} \leq \| \vert D \vert^{n\alpha} Q\|_{L^2} + \| \vert D \vert^{n\alpha}(Q^p) \|_{L^2} \leq C \left( 1 + \| Q\|_{L^\infty}^{p-1}\right) \| Q \|_{H^{n\alpha}} <\infty. 
\end{align*}

\begin{rema}
Instead of the set of assumptions \eqref{hypo:Q} on $Q$, one can ask $Q\in L^2(\R) \cap L^{p+1}(\R)$ to obtain the same result (see Lemma B.1 of \cite{FL13}).
\end{rema}

\textbf{Acknowledgements} The authors thank Martin Oen Paulsen for interesting discussions and Didier Pilod for his constant support on this project. The authors were supported by a Trond Mohn foundation grant.


\bibliographystyle{plain}
\bibliography{Z_biblio}

\end{document}